\documentclass[12pt]{amsart}

\title{A Kurosh-Type Theorem For Type III Factors}
\author{Jason Asher}

\usepackage{amsmath}
\usepackage{amssymb}
\usepackage{fullpage}

\newtheorem{thm}{Theorem}[section]
\newtheorem{lem}[thm]{Lemma}
\newtheorem{prop}[thm]{Proposition}
\newtheorem{cor}[thm]{Corollary}

\newcommand{\nolineover}[2]{\genfrac{}{}{0pt}{}{#1}{#2}}
\makeatletter
\def\freeprodsize@{\@setfontsize\freeprodsize\@xxpt\@xxpt}
\def\freeprod@{\mathop{\hbox{\freeprodsize@ $*$}}}
\def\freeprod{\freeprod@\displaylimits}
\makeatother
\newcommand{\bigast}{\freeprod}

\begin{document}
\maketitle

\section{Introduction}

In \cite{ozkurosh}, Ozawa obtained analogues of the Kurosh subgroup theorem (and its consequences) in the setting of free products of semiexact $\text{II}_1$ factors with respect to the canonical tracial states.  In particular he was able to prove a certain unique-factorization theorem that distinquishes, for example, the n-various $L(F_\infty) \ast (L(F_\infty) \otimes R)^{\ast n}$. The paper was a continuation of the joint work \cite{opprime} with Popa that proved various unique-factorization results for tensor products of $\text{II}_1$ factors.

These papers shared a particular combination of ideas from \cite{ozsolid} and \cite{postrong}. The first was a $C^*$-algebraic method for detecting injectivity of a von Neumann algebra, adopted from its original context of proving solidity of finite von Neumann algebras with the Akemann-Ostrand property. This method was used in concert with Popa's intertwining-by-bimodules technique that (roughly speaking) allows one to conclude unitary conjugacy results from a spatial condition. More precisely, if $A$ and $B$ are diffuse subalgebras of a finite von Neumann algebra $M$, the presence of an $A$-$B$ sub-bimodule of $L^2(M)$ with finite right $B$-dimension implies that a `corner' of $A$ can be conjugated into $B$ by a partial isometry in $M$. The technique was actually shown by Popa to apply to the case of $(M,\varphi)$ not necessarily finite (but with discrete decomposition), though one needs the additional assumption that $A$ and $B$ are subalgebras of the centralizer $M^\varphi$.

In this paper we will demonstrate how the above results can be extended to the case of free products of semiexact $\text{II}_1$ factors with respect to arbitrary (non-tracial) states. The primary difficulty this setting presents is that the resulting factors will be necessarily of type $\text{III}$. Following the same basic outline of the proof from \cite{ozkurosh}, we will therefore utilize a generalization of Popa's intertwining technique that allows the ambient algebra $M$ to be completely arbitrary, and makes no assumption about the relative position of the subalgebras $A$ and $B$ in $M$. One of the fruits of our labors will be the result that for $M_i$ semiexact non-prime non-injective $II_1$ factors with faithful normal states $\varphi_i$, the reduced free product $\bigast_{i=1}^m (M_i , \varphi_i)$ of the $(M_i,\varphi_i)$ is uniquely written in such a way as the free product of such factors, up to permutation of indices and stable isomorphism.

Acknowledgement: The author would like to thank Dimitri Shlyakhtenko for initially suggesting this investigation, and for his invaluable subsequent advice and direction.

\section{Reduced Free Product}

We will first review the reduced free product construction \cite{vo} in our particular case, retaining the notation from \cite{ozkurosh}. For $i \in \{1,\dots,m\}$, let $M_i$ be semiexact $\text{II}_1$ factors with faithful normal states $\varphi_i$, and let $H_i = L^2(M_i,\varphi_i)$. Denote by $\xi_i$ the $M_i$-cyclic-and-separating state vector $\hat 1 \in H_i$, and let $J_i$ be the modular conjugation on $H_i$. Finally, let $H_i^0 = H_i \ominus \mathbb{C}\xi_i$, and construct a new Hilbert space $H$ as follows:
$$H = \mathbb{C}\xi \oplus \bigoplus_{n \ge 1}\bigoplus_{\nolineover{i_1,\dots,i_n \in \{1,\dots,m\}} {i_1 \ne i_2 \ne \dots \ne i_n}} H_{i_1}^0 \otimes \cdots \otimes H_{i_n}^0.$$

Next, let $H(i)$ be the subspace of $H$ given by the direct sum of $\mathbb{C}\xi$ and those terms in the direct sum above with $i_1 \ne i$. For each $i$ there is a unitary isomorphism
$$U_i : H \rightarrow H_i \otimes H(i)$$
which identifies $H(i)$ with $\mathbb{C}\xi \otimes H(i)$ and $H(i)^\perp$ with $H_i^0 \otimes H(i)$. Define a representation $\lambda_i:\mathbb{B}(H_i) \rightarrow \mathbb{B}(H)$ by
$$\lambda_i(a) = U_i^*(a \otimes 1_{H(i)}) U_i.$$
The reduced free-product $(M,\varphi) = \bigast_{i=1}^m (M_i , \varphi_i)$ is then the von Neumann algebra in $\mathbb{B}(H)$ generated by $\lambda_i(M_i)$ for $i \in \{1,\dots,m\}$, with faithful normal state $\varphi$ given by
$$\varphi(x) = \left< x\xi,\xi \right>_H.$$

Note that $\xi$ is  an $M$-cyclic-and-separating state vector, and so $H = L^2(M,\varphi)$ with the associated modular conjugation $J$ given by
$$J(\eta_{i_1} \otimes \cdots \otimes \eta_{i_n}) = J_{i_n}\eta_{i_n} \otimes \cdots \otimes J_{i_1}\eta_{i_1}$$
where $\eta_{i_1} \otimes \cdots \otimes \eta_{i_n} \in H_{i_1}^0 \otimes \cdots \otimes H_{i_n}^0$ as in the direct sum decomposition above. Furthermore, one can thus see that $JH(i)$ is the direct sum of $\mathbb{C}\xi$ and those terms in the direct sum above with $i_n \ne i$.

For each $i$ there is thus another unitary isomorphism 
$$V_i:H \rightarrow JH(i) \otimes H_i$$ 
which identifies $JH(i)$ with $JH(i) \otimes \mathbb{C}\xi$ and $JH(i)^\perp$ with $JH(i) \otimes H_i^0$. Then, precisely as in \cite{ozkurosh} Lemma 3.1, we have the following 
\begin{lem}
For $a \in \mathbb{B}(H_i)$, $j\ne i$, and with $P_\xi$ the projection of $H$ onto $\mathbb{C}\xi$ we have
\begin{align*}
\lambda_i(a) &= JV_i^*(1_{JH(i)} \otimes J_i a J_i)V_i J \\
&= V_i^*(P_\xi \otimes a + \lambda_i(a)_{|_{JH(i)\ominus \mathbb{C}\xi}} \otimes 1_{H_i})V_i \\
&= V_j^*(\lambda_i(a)_{|_{JH(i)}} \otimes 1_{H_i})V_j
\end{align*}
\end{lem}

This shows that $V_i$ is a right $M_i$-module isomorphism of $L^2(M,\varphi)$ with $K_i\otimes H_i$, where $K_i = JH(i)$. Furthermore, if we let $D_i = V_i^*(\mathbb{K}(K_i) \otimes_{min} \mathbb{B}(H_i))V_i$ then this lemma also shows that $\lambda_j(a)D_i \subseteq D_i$ for $i,j \in \{1,\dots,m\}$.

Let $B_i \subseteq M_i$ be a $\sigma$-weakly dense exact $C^*$ algebra, and let $B$ be the $C^*$ algebra in $\mathbb{B}(H)$ generated by $\lambda_i(B_i)$ for $i \in \{1,\dots,m\}$. Then for $C = JBJ$ we have the following lemma, as in \cite{ozkurosh}, Prop. 3.2.

\begin{lem}
Suppose that $\Psi: \mathbb{B}(H) \rightarrow \mathbb{B}(H)$ is a $C$-bimodule unital completely positive map. If $D_i \subset ker \Psi$ for $i\in \{1,\dots,m\}$, then the unital completely positive map
\begin{align*}
\tilde\Psi : B \otimes C \, &\rightarrow \ \mathbb{B}(H) \\
\sum_{k=1}^n a_k \otimes x_k &\mapsto \Psi\left(\sum_{k=1}^n a_k x_k\right)\\
\end{align*}
is continuous with respect to the minimal tensor norm.
\end{lem}

Finally, we will need the following control of partial isometries that normalize diffuse subalgebras in free products, see \cite{poorth} and \cite{ozkurosh} Lemma 2.3. For completeness, we include the proof.

\begin{lem}
Let $(M,\varphi) =  \bigast_{i=1}^l (M_i , \varphi_i)$ be the reduced free product of $\text{II}_1$ factors $M_i$ with faithful normal states $\varphi_i$, and let $Q \subseteq M_j$ be a diffuse von Neumann subalgebra. If $v \in M$ is a nonzero partial isometry with $vv^* \in Q' \cap M$ and such that $v^* Q v \subseteq M_i$, then $i=j$ and $v \in M_j$.
\end{lem}

Proof: Let $\zeta \in L^2(M_i,\varphi_i) \subseteq L^2(M,\varphi)$ be such that $\left < \cdot \zeta, \zeta \right>$ gives the canonical trace on $M_i$. Note that then $x\zeta = J_\varphi x^* J_\varphi \zeta = \zeta x$ for all $x \in M_i$. If we let $\eta = v\zeta$, then $\eta \ne 0$ and
$$q \eta = vv^* q v \zeta = v \zeta (v^* q v) = \eta (v^* q v)$$
for all $q \in Q$.

Suppose, by contradiction, that $i \ne j$. Then $L^2(M,\varphi)$ as an $M_j - M_i$ bimodule is an infinite multiple of the `coarse' bimodule $L^2(M_j,\varphi_j) \otimes L^2(M_i,\varphi_i)$. Considering the projection of $\eta$ onto any coarse component, we get a Hilbert-Schmidt operator $T$ on $L^2(M,\varphi)$ such that
$$q T = T (v^*qv)$$
for all $q \in Q$. As $Q$ is diffuse, we must have that $Tv^* = 0$, so that $T = Tv^*v = 0$. Thus, $\eta = 0$, a contradiction.

So, we must have that $i = j$. Now, note that $L^2(M,\varphi)$ as an $M_j - M_j$ bimodule is isomorphic to the direct sum of $L^2(M_j,\varphi_j)$ and an infinite multiple of the coarse bimodule $L^2(M_j,\varphi_j) \otimes L^2(M_j,\varphi_j)$. By the above argument, we must have that $\eta \in L^2(M_j,\varphi_j)$. Thus, considering these $L^2$-vectors as unbounded operators affiliated with $M_j$, we get that $v = \eta\zeta^{-1} \in M_j$. $\square$

\section{Intertwining-by-bimodules}
We present here a generalization of the intertwining-by-bimodules technique of Popa, see Theorem 2.1 of \cite{postrong}. The technical advantage is that we allow the ambient algebra to be arbitary, and we make no assumption about the location of the subalgebras.

The following lemma is standard - for a reference on this and other aspects of the modular theory of Tomita and Takesaki see \cite{tavol2}.
\begin{lem}
Let $M$ be a von Neumann algebra in a standard representation $(M,H,J,P)$, where $J$ is the modular conjugation, and $P$ is the self-dual cone. Then any vector $\eta \in H$ has a unique polar decomposition as $\eta = u |\eta|$ for $|\eta| \in P$ and a partial isometry $u \in M$ such that the source of $u$ is $[M' |\eta|]$, and the target of $u$ is $[M' \eta]$. Furthermore if $\eta$ satisfies an equation of the form $x \eta = J x^* J \eta$ for all $x$ in a (possibly nonunital) subalgebra $A$ of $M$, then $|\eta|$ satisfies the same equation, and $u$ satisfies $xu=ux$ for all $x \in A$.
\end{lem}
Proof:
Note that $(M',H,J,P)$ is a standard representation for $M'$, and consider the positive normal functional $\varphi_\eta$ on $M'$ given by
$$\varphi_\eta(z) = \left< z \eta, \eta \right>.$$
Then there is a unique vector $|\eta| \in B$ such that $\varphi_\eta(z) = \left< z |\eta|, |\eta| \right>$. Defining $u$ on $M' |\eta|$ by $u z |\eta| = z \eta$ is well-defined (and when extended by continuity gives a partial isometry) because
$||z \eta||^2 = \varphi_\eta(z^*z) = ||z |\eta|||^2$. Then $u$ is clearly in $M$ and furthermore, uniqueness is now obvious.

Suppose further that $\eta$ satisfies $x \eta = J x^* J \eta$ for all $x \in A$. By considering $A + \mathbb{C}(1-p)$ for $p$ the unit of $A$, we may assume that $A$ is unital. Then for any unitary $w \in A$ one has that
$$(wuw^*) (wJwJ |\eta|) = wJwJ \eta = \eta.$$
Also, $\zeta = wJwJ |\eta| \in B$, and $wuw^*$ is a partial isometry with source $[M' \zeta]$ and target $[M' \eta]$. Thus, by uniqueness, we obtain the result.

$\square$

\begin{prop}
Let $M$ be a von Neumann algebra with faithful normal state $\varphi$, and let $A, B \subseteq M$ be two factors of type $II_1$. Suppose that there is a nonzero $A$-$B$ bimodule $H \subset L^2(M,\varphi)$ such that the von Neumann dimension of H as a right B module satisfies $\dim_B H_B < \infty$. Then there are nonzero projections $q \in A$ and $p \in B$, an injective *-homomorphism $\theta: qAq \rightarrow pBp$, and a nonzero partial isometry $v \in M$ such that $v^*v \le p$, $vv^* \le q$, and $xv = v\theta(x)$ for all $x \in qAq$.
\end{prop}

Proof:

We can pick nonzero projections $q \in A$ and $p \in B$ such that $K = qJ_\varphi p J_\varphi H$ is a nonzero $qAq$-$pBp$ bimodule with $\dim_{pBp} K_{pBp} = 1$. Then, we can find an right Hilbert $pBp$-module isomorphism $\alpha: L^2(pBp,\tau_{pBp}) \rightarrow H$, and
define $\theta: qAq\rightarrow pBp$ by
$$\theta(x) \zeta = \alpha^{-1}(x\alpha(\zeta) ) \text{ for all }\zeta \in L^2(pBp,\tau_{pBp}).$$
Note that $\theta$ is a nonzero, and hence injective, *-homomorphism.

Furthermore, if we let $\xi = \alpha(\hat p)$, then for $x \in qAq$,
$$x \xi = \alpha(\theta(x) \hat p)
= \alpha(\hat p \theta(x))
= J_\varphi \theta(x)^* J_\varphi \xi.$$

Consider now $M_2 \otimes M$ with the faithful normal state $\tilde\varphi = \operatorname{Tr} \otimes \varphi$, and define for $x \in qAq$,
$$ \tilde x =
\left[
\begin{array}{cc}
x & 0\\
0 & \theta(x)
\end{array}
\right].
$$
Also, let $\tilde\xi \in L^2(M_2 \otimes M,\tilde\varphi)$ be defined as
$$ \tilde\xi =
\left[
\begin{array}{cc}
0 & \xi\\
0 & 0
\end{array}
\right]
$$
and note that we then have the equation
$$\tilde x \tilde \xi = J_{\tilde\varphi} \tilde x^* J_{\tilde\varphi} \tilde\xi$$
for all $x \in qAq$.

By the lemma above, if we let $\tilde \xi = u |\tilde \xi|$ be the polar decomposition of $\tilde \xi$ then $u$ in $M_2 \otimes M$ is a partial isometry satisfying
$$ \tilde x u = u \tilde x$$
for all $x \in qAq$.
It is easy to see that $u$ has the form
$$ u=
\left[
\begin{array}{cc}
0 & w\\
0 & 0
\end{array}
\right]
$$
for a partial isometry $w \in M$ that satisfies
$xw = w\theta(x)$ for all $x \in qAq$. In fact $w$ is the polar part of $\xi$, and so letting $v = qw = qw\theta(q) = w\theta(q)$ , we obtain the desired result.

$\square$

\section{Kurosh-Type Theorem}
We will now present the main technical result, a generalization of Theorem 3.3 in \cite{ozkurosh}. Note that much of the argument remains the same, but we include it for narrative coherence. Also, in everything that follows, all von Neumann algebras are assumed to be separable.

\begin{thm}
Let $M_i$, for $i\in \{1,\dots,l\}$, be semiexact $\text{II}_1$ factors with faithful normal states $\varphi_i$, and let $(M,\varphi) =  \bigast_{i=1}^l (M_i , \varphi_i)$ be their reduced free product. Suppose that $Q \subseteq M$ is an injective $\text{II}_1$ factor with a normal conditional expectation $E_Q: M \rightarrow Q$. If $Q' \cap M$ is non-injective, then there exists an index $i$ and a maximal partial isometry $u\in M$ with $uu^* \in Q' \cap M$ and such that $u^*Qu \subseteq M_i$.
\end{thm}

Proof:
Recall all of the notation from the reduced free product construction above.
Let 
$$\Psi_Q: \mathbb{B} (L^2(M,\varphi)) \rightarrow Q' \cap \mathbb{B}(L^2(M,\varphi))$$ 
be a proper conditional expectation, i.e. such that for all $x \in \mathbb{B}(L^2(M,\varphi))$ we have 
$$\Psi_Q(x) \in \overline{\text{co}}^{\omega} \{ uxu^* : u\in U(Q)\}.$$ Note that $\Psi_Q$ is a $Q' \cap \mathbb{B}(L^2(M,\varphi))$-bimodule map.

As $Q$ is in the range of a normal conditional expectation from $M$, we can find a faithful normal state $\nu$ on $M$ which centralizes $Q$. Then ${\Psi_Q}_{|_M}$ is $\nu$-preserving, and so ${\Psi_Q}_{|_M}$ is the unique $\nu$-preserving normal conditional expectation $E_{Q' \cap M} : M \rightarrow Q'\cap M$.
Lemma 5 from \cite{ozsolid} thus implies that, as $Q' \cap M$ is non-injective, the map
\begin{align*}
\tilde\Psi_Q : B \otimes C \, &\rightarrow \ \mathbb{B}(L^2(M,\varphi)) \\
\sum_{k=1}^n a_k \otimes x_k &\mapsto \Psi_Q\left(\sum_{k=1}^n a_k x_k\right)\\
\end{align*}
is not continuous with respect to the minimal tensor norm. In turn, Lemma 2.2 above implies that there must be an index $i$ such that $D_i \not\subset \text{ker}(\Psi_Q)$.

By continuity of $\Psi_Q$, there must be a finite rank projection $f \in \mathbb{B}(K_i)$ such that $x = \Psi_Q(f \otimes 1) \ne 0$, where we have identified $L^2(M,\varphi)$ with $K_i \otimes H_i$ via $V_i$. Note that $x$ commutes with the right $M_i$-action, and so $x\in (\mathbb{B}(K_i) \bar\otimes M_i) \cap Q'$. Furthermore, $\text{Tr} \otimes \tau_{M_i}(x) \le \text{Tr}(f) < \infty$, and so by taking a suitable spectral projection $e$ of $x$ we get that $eL^2(M,\varphi)$ is a nonzero $Q$-$M_i$ bimodule with $\dim_{M_i} eL^2(M,\varphi) < \infty$.

Thus, by Proposition 3.2, there are nonzero projections $q \in Q$ and $p \in M_i$, a *-homomorphism $\theta: qQq \rightarrow pM_ip$, and a nonzero partial isometry $v \in M$ such that $v^*v \le p$, $vv^* \le q$ and for all $a \in qQq$ we have
$$av = v\theta(a).$$
Note that this implies that $vv^* \in (qQq)' \cap qMq = q(Q' \cap M)q$ and so we can write $vv^* = qq'$ for $q' \in Q' \cap M$. Furthermore, this also implies that $v^*v \in \theta(qQq)' \cap pMp$, and so by Lemma 2.3 we have that $v^*v \in pM_i p$. Thus for all $a \in Q$,
$$v^* a v = v^*v \theta(qaq) \in M_i.$$

Also, note that all of this remains true after trimming $q$ and $q'$ (and $v$, $\theta$, in turn). So, by restricting $q$ (and hence $v$ and $\theta$) we may assume that $\tau_Q(q) = \frac{1}{n}$ for some integer $n$. Then, by restricting $q'$ (and hence $v$), we may assume that $\tau_{M_i}(v^*v) = \frac{1}{mn}$ for some integer $m$. We can then find $n$ partial isometries $u_i \in Q$ and $w_i \in M_i$, respectively, such that $\sum_i u_i u_i^* = 1$, $u_i^*u_i = q$, $\sum_i w_i^* w_i \le 1$, and $w_i w_i^* = v^*v$. Then, letting $w = \sum_i u_i v w_i$ we get that for $a \in Q$, 
$$w^* a w = \sum_i \sum_j w_i^* v^* u_i^* a u_j v w_j \in M_i$$
while
$$w w^* = \sum_i u_i q' u_i^* = q'  \in Q' \cap M.$$

Note that $\tau_{M_i}(w^*w) = \frac{1}{m}$. If $Q' \cap M$ is infinite, then by repeating an argument like the one above, we can enlarge $w$ to a partial isometry $u$ with $u^*u = 1$, $uu^* \in Q' \cap M$ and such that $u^*Qu \subseteq M_i$. Furthermore, if $Q' \cap M$ is type $\text{III}$ (or if $q'$ is an infinite projection and $Q' \cap M$ is type $\text{II}_\infty$) we can easily dilate $u$ to a unitary.

Otherwise, if $Q' \cap M$ is type $\text{II}_1$ we have two cases. First, if $\tau_{Q'\cap M}(q') \le \frac{1}{m}$, then we can proceed precisely as in the infinite situation and obtain $u$ with the same properties. If on the other hand $\tau_{Q'\cap M}(q') > \frac{1}{m}$, we can restrict $w$ and assume that $\tau_{Q'\cap M}(ww^*) = \frac{1}{m} \ge \tau_{M_i}(w^*w)$. Then, proceeding as above we can enlarge $w$ to a partial isometry $u$ with $uu^*=1$ and such that $u^*Qu \subseteq M_i$. $\square$

The following lemma follows easily from some well-known results of Popa, see A.1.1 and A.1.2 of \cite{poclass}. The proof in the case that the ambient algebra $M$ is finite can be found as Proposition 13 of \cite{opprime}, and the general case follows by the exact same argument where one replaces trace-preserving conditional expectations with those that preserve a state $\varphi$ on $M$ that centralizes $N$.

\begin{lem}
Let $N \subseteq M$ be an inclusion of factors with $N$ type $\text{II}_1$ and in the range of a normal conditional expectation $E_N:M \rightarrow N$. Then there is an injective $\text{II}_1$ factor $Q \subseteq N$ such that $Q' \cap M = N' \cap M$.
\end{lem}

\begin{cor}
Let $M_i$, for $i\in \{1,\dots,l\}$, be semiexact $\text{II}_1$ factors with faithful normal states $\varphi_i$, and let $(M,\varphi) =  \bigast_{i=1}^l (M_i , \varphi_i)$ be their reduced free product. Suppose that $N \subseteq M$ is an non-prime non-injective $\text{II}_1$ factor with a normal conditional expectation $E_N: M \rightarrow N$. Then there is an index $i$ and a  maximal partial isometry $u$ such that $u^*Nu \subseteq M_i$.
\end{cor}
Proof: Write $N = N_1 \otimes N_2$ for $\text{II}_1$ factors $N_i$ with $N_2$ non-injective. By Lemma 4.2, there is an injective $\text{II}_1$ factor $Q \subseteq N_1$ such that $Q' \cap M$ = $N_1' \cap M$. Note that $N_2 \subseteq Q' \cap M$, and thus $Q' \cap M$ is non-injective. Furthermore, $Q$ is in the range of a normal conditional expectation from $M$. Thus, by Theorem 4.1, there is an index $i$ and a maximal partial isometry $u$ such that $uu^* \in Q' \cap M = N_1' \cap M$, and $u^*Qu \subseteq M_i$.

Now, note that if $x\in N_2$ and $q \in Q$, then
$$(u^* x u)( u^* q u) = u^* x q u = u^* q x u = (u^* q u)( u^* x u)$$
and so by Lemma 2.3, $u^* x u \in M_i$.

Furthermore, if $y\in N_1$ and $x \in N_2$, then as $uu^* \in N_1' \cap M$,
$$(u^* y u)( u^* x u) = u^* y x u = u^* x y u = (u^* x u)( u^* y u)$$
and so by Lemma 2.3 again, $u^* y u \in M_i$.

Thus, $u^* N_1 N_2 u = u^*N_1u u^*N_2 u \subseteq M_i$, and so $u^*Nu \subseteq M_i$. $\square$

\begin{cor}
Let $\{M_i\}_{i=1}^m$, $\{N_j\}_{j=1}^n$ be semiexact non-prime non-injective $\text{II}_1$ factors with faithful normal states $\varphi_i$, $\psi_j$, and let $M,\varphi$, $N, \psi$ be their respective reduced free products. If $M \simeq N$, then $m = n$ and after reordering indices $M_i$ is isomorphic to an amplification of $N_i$.
\end{cor}
Proof: Identify $M$ and $N$ via the implied isomorphism. By the proof of Corollary 4.3, for each index $i$ there is an injective $\text{II}_1$ factor $Q_i \subseteq M_i$, an index $\alpha(i)$, and a maximal partial isometry $u_i \in M$ such that $u_iu_i^* \in Q_i' \cap M$ and $u_i^* M_i u_i \subseteq N_{\alpha(i)}$. Similarly for each index $j$ there is an injective $\text{II}_1$ factor $P_j \subseteq N_j$, an index $\beta(j)$, and a maximal partial isometry $v_j \in M$ such that $v_jv_j^* \in P_j' \cap M$ and $v_j^*N_j v_j \subseteq M_{\beta(j)}$. 

Consider the projections $p_i = u_i^*u_i$ and $q_i = v_{\alpha(i)}{v_{\alpha(i)}}^*$ in $N_{\alpha(i)}$ (Lemma 2.3 implies that $q_i \in N_{\alpha(i)}$). First, by restricting the target projection of $u_i$ in $Q_i' \cap M$, we may assume that $\tau_{N_{\alpha(i)}}(p_i) \le \tau_{N_{\alpha(i)}}(q_i)$. Furthermore, by factoriality of $N_{\alpha(i)}$ we can rotate $v_i$ by a partial isometry to assure that $p_i \le q_i$. Then, if we let $w_i = u_i v_{\alpha(i)}$, $w_i$ is a partial isometry such that $w_i^* Q_i w_i \subseteq M_{\beta(\alpha(i))}$, and $w_iw_i^* = u_i u_i^* \in Q'_i \cap M$. Thus, by Lemma 2.3, $\beta(\alpha(i)) = i$ and $w_i \in M_i$.
So, $u_i$ implements an isomorphism of $(u_i u_i^*)A(u_i u_i^*)$ with $pN_\alpha(i)p$.

A completely symmetric argument shows that $\alpha(\beta(j))=j$. $\square$

Note that this last corollary implies that, for example, if $d_i \in M = L(\mathbb{F}_n) \otimes L(\mathbb{F}_m)$ are positive and $\varphi_i(\cdot) \sim \tau_M(d_i \cdot)$, then $\bigast_{i=1}^l (M, \varphi_i)$ are non-isomorphic for different $l$.


\begin{thebibliography}{article}

%
\bibitem[Oz]{ozkurosh}
N. Ozawa,
\textit{A Kurosh type theorem for type $\mathrm{II}_1$ factors.}
Int. Math. Res. Not. (2006), doi:10.1155/IMRN/2006/97560
%
\bibitem[Oz2]{ozsolid}
N. Ozawa, 
\textit{Solid von Neumann algebras.} 
Acta Math. \textbf{192} (2004), 111--117.
%
\bibitem[OP]{opprime}
N. Ozawa and S. Popa, 
\textit{Some prime factorization results 
for type~$\mathrm{II}_1$ factors.}
Invent. Math. \textbf{156} (2004), 223--234.
%
\bibitem[Po]{poclass}
S. Popa,
\textit{Classification of amenable subfactors of type $\mathrm{II}$.}
Acta Math. \textbf{172} (1994) 163--255.
%
\bibitem[Po2]{poorth}
S. Popa,
\textit{Orthogonal pairs of *-subalgebras in finite von Neumann algebras.} 
J. Operator Theory \textbf{9} (1983),  253--268.
%
\bibitem[Po3]{postrong}
S. Popa, 
\textit{Strong rigidity of $\mathrm{II}_1$ factors coming 
from malleable actions of weakly rigid groups, I.}
Invent. Math. \textbf{165} (2006), 369--408. 
%
\bibitem[Ta]{tavol2}
M. Takesaki,
\textit{Theory of Operator Algebras $\mathrm{II}$.} 
Encyclopedia of Mathematical Sciences \textbf{125}, Springer--Verlag, 2000. 
%
\bibitem[Vo]{vo}
D. Voiculescu,
\textit{Symmetries of some reduced free product $C^*$--algebras.}
Operator Algebras and Their Connections with Topology and Ergodic Theory, Lecture Notes in Mathematics, vol. 1132, Springer--Verlag, 1985, 556--588.
%

\end{thebibliography}
\end{document}